\documentclass[11pt]{amsart}
\usepackage{graphicx}
\usepackage[centertags]{amsmath}
\usepackage{amsfonts}
\usepackage{amssymb}
\usepackage{amsthm}
\usepackage{newlfont}
\newtheorem{theorem}{Theorem}[section]

\newtheorem{prop}[theorem]{Proposition}
\theoremstyle{definition}

\theoremstyle{remark}

\newcommand{\tr}{\textrm{tr}}
\newcommand{\id}{\textrm{id}}

\title{Curvature properties of $4$-dimensional\\ Riemannian manifolds with a circulant\\ structure}

\author{Iva Dokuzova}

\begin{document}

%

%%%%%%%%%%%%%%%%%%%%%%%%%%%%%%%%%
\begin{abstract}
%%%%%%%%%%%%%%%%%%%%%%%%%%%%%%%%%%%
We consider a $4$-dimensional Riemannian manifold $M$ equip\-ped with a circulant structure $q$, which is an isometry with respect to the metric $g$ and $q^{4}=\id$, $q^{2}\neq \pm \id$. For such a manifold $(M, g, q)$ we obtain some assertions for the sectional curvatures of $2$-planes. We construct an example of such a manifold on a Lie group and we find some of its geometric characteristics.
\end{abstract}

\maketitle
%%% ----------------------------------------------------------------------


%%%%%%%%%%%%%%%%%%%%%%%%%%%%%%%%%%%%%%%%%%%%%% 2
\section*{Introduction}
The circulant matrices are well-studied (for example \cite{8}, \cite{152}). They have application to Vibration analysis, Graph theory, Linear codes, Geometry (for example \cite{13}, \cite{72}, \cite{2}).

  The study of manifolds with additional structures plays an important role in differential geometry. In such manifolds substantial results are associated with the sectional curvatures of some characteristic $2$-planes of the tangent space of the manifolds (for example \cite{22}, \cite{132}, \cite{52}).

 In the present paper we consider some curvature properties of $4$-dimen\-sional Riemannian manifolds with a circulant structure $q$ with $q^{4}=\id$, which is an isometry with respect to the metric $g$. We continue research made in \cite{1} for such manifolds and construct an example of these manifolds.

The paper is organized as follows. In Sect. \ref{sec:1} we give some necessary facts from \cite{1} about a $4$-dimensional differentiable manifold $M$ with a Riemannian metric $g$, equipped with a circulant structure $q$, which is an isometry with respect to the metric $g$ and $q^{4}= \id$, $q^{2}\neq \pm \id$.
In Sect.~\ref{sec:2} we establish that the sectional curvatures of the $2$-planes $\{u, qu\}$ and $\{u, q^{2}u\}$ are expressed by the angles $\angle(u, qu)$ and $\angle(u, q^{2}u)$, respectively. The main results here are Theorem~\ref{th6} and Theorem~\ref{th5}. We obtain relations between the sectional curvatures of some characteristic $2$-planes in the tangent space on the manifold $(M, g, q)$.
In Sect.~\ref{sec:3} we construct an example of such a manifold on a Lie group and we find some of its geometric characteristics.

\section{Preliminaries}\label{sec:1}

In this section we recall facts from \cite{1}, which are necessary for our future consideration.

Let $M$ be a $4$-dimensional Riemannian manifold with a metric $g$. Let $q$ be an endomorphism in the tangent space $T_{p}M$, $p\in M$ on the manifold $M$ with local coordinates given by the circulant matrix
\begin{equation}\label{f4}
    (q_{i}^{j})=\begin{pmatrix}
      0 & 1 & 0 & 0\\
      0 & 0 & 1 & 0 \\
      0 & 0 & 0 & 1\\
      1 & 0 & 0 & 0\\
    \end{pmatrix}.
\end{equation}
Then
\begin{equation}\label{q4}
    q^{4}=\id,\qquad q^{2}\neq\pm \id.
\end{equation}
We suppose that $g$ is positive definite metric and the structure $q$ of the manifold $M$ is an isometry with respect to the metric $g$, i.e.
\begin{equation}\label{2.1}
    g(qx, qy)=g(x, y).
\end{equation}
Anywhere in this work $x, y, z, u$ will stand for arbitrary elements of the algebra of the smooth vector fields on $M$ or vectors in the tangent space $T_{p}M$. The Einstein summation convention is used, the range of the summation indices being always $\{1, 2, 3, 4\}$.

We denote by $(M, g, q)$ the manifold $M$ equipped with the metric $g$ and the structure $q$.

Easily finding that \eqref{f4} and \eqref{2.1} imply a circulant matrix of components of $g$.

%%% ----------------------------------------------------------------

A basis of type $\{x, qx, q^{2}x, q^{3}x\}$ of $T_{p}M$ is called a $q$-\textit{basis}. In this case we say that \textit{the vector $x$ induces a $q$-basis of} $T_{p}M$.

If a vector $x$ induces a $q$-basis, then for the angles
 $\angle(x,qx),\ \angle(x,q^{2}x)$, $\angle(qx,q^{2}x)$, $\angle(qx,q^{3}x)$, $\angle(x,q^{3}x)$ and $\angle(q^{2}x,q^{3}x)$
 we have
 \begin{equation*}
    \angle(x,qx)=\angle(qx,q^{2}x)=\angle(x,q^{3}x)=\angle(q^{2}x,q^{3}x),\quad  \angle(x,q^{2}x)=\angle(qx,q^{3}x).
\end{equation*}

In our further research we will use an orthogonal $q$-basis. The existence of such bases is proved in \cite{1}.

\section{Some curvature properties}\label{sec:2}

Let $\nabla$ be the Riemannian connection of the metric $g$ on $(M, g, q)$. The curvature tensor $R$ of $\nabla$ is determined by $R(x, y)z=\nabla_{x}\nabla_{y}z-\nabla_{y}\nabla_{x}z-\nabla_{[x,y]}z$. The tensor of type $(0, 4)$  associated with $R$ is defined as follows
$$R(x, y, z, u)=g(R(x, y)z,u).$$

If we denote $P=q^{2}$, then the conditions \eqref{q4} and \eqref{2.1} imply $P^{2}=\id,$ $\ P\neq \pm \id$, $g(Px, Py)=g(x,y)$.
Thus, $(M, g, P)$ is a Riemannian manifold with an almost product structure $P$. It follows from \eqref{f4} that $\tr P=0$. For such manifolds is valid Staikova-Gribachev classification (\cite{51}). The class  $W_{0}$ defined by $\nabla P=0$ in this classification is common to all classes. Every manifold in this class satisfies the identity $R(x, y, Pz, Pu)=R(x, y, z, u)$.
In \cite{1} it is proved analogous identity
\begin{equation}\label{R1}
  R(x, y, qz, qu)=R(x, y, z, u),
\end{equation}
for a manifold $(M, g, q)$ with the condition $\nabla q=0$.
By using \eqref{R1} and the symmetries of $R$, it is easy to find that
\begin{equation}\label{R}
  R(qx, qy, qz, qu)=R(x, y, z, u).
\end{equation}
Since the latter equality follows from \eqref{R1}, the class of manifolds $(M, g, q)$ with the condition \eqref{R} is more general than the class $(M, g, q)$ with the condition \eqref{R1}.

If $\{x, y\}$ is a non-degenerate $2$-plane spanned by vectors $x, y \in T_{p}M$, then its sectional curvature is (\cite{11})
\begin{equation}\label{3.3}
    \mu(x,y)=\frac{R(x, y, x, y)}{g(x, x)g(y, y)-g^{2}(x, y)}\ .
\end{equation}
\begin{theorem}\label{predl}
Let $(M, g, q)$ be a manifold with property \eqref{R}. If a vector $x$ induces a $q$-basis, then for the sectional curvatures of the basic $2$-planes we have
\begin{equation}\label{3.21}
\mu(x,qx)=\mu(qx,q^{2}x)=\mu(q^{2}x,q^{3}x)=\mu(q^{3}x,x),
\end{equation}
\begin{equation}\label{3.31}
\mu(x,q^{2}x)=\mu(qx,q^{3}x).
\end{equation}
\end{theorem}
\begin{proof}
From (\ref{R}) we have
\begin{equation}\label{Rv1}
R(x,y,z,u)=R(qx,qy,qz,qu)=R(q^{2}x,q^{2}y,q^{2}z,q^{2}u).
\end{equation}
In \eqref{Rv1} we substitute

1) $qx$ for $y$, $x$ for $z$, $qx$ for $u$, and we get
 \begin{equation}\label{slRv1}
 R(x,qx,x,qx)=R(qx,q^{2}x,qx,q^{2}x )=R(q^{2}x,q^{3}x,q^{2}x,q^{3}x),
 \end{equation}

2) $q^{3}x$ for $y$, $x$ for $z$, $q^{3}x$ for $u$, and we obtain
  \begin{equation}\label{dop2}
    R(x, q^{3}x, x, q^{3}x)=R(x, qx, x, qx),
    \end{equation}

3) $q^{2}x$ for $y$, $x$ for $z$, $q^{2}x$ for $u$, then
 \begin{equation}\label{dop4}
 R(x,q^{2}x,x,q^{2}x)=R(qx,q^{3}x,qx,q^{3}x).
 \end{equation}
 The equality \eqref{3.21} follows from \eqref{2.1}, \eqref{3.3}, \eqref{slRv1} and \eqref{dop2}. In a similar way, from \eqref{2.1}, \eqref{3.3} and \eqref{dop4} we get \eqref{3.31}.
\end{proof}

Let $x$ induce a $q$-basis $\{x, qx, q^{2}x, q^{3}x\}$. Due to Theorem~\ref{predl} there are only two different basic sectional curvatures. Therefore, we consider only the sectional curvatures $\mu(x, qx)$ and $\mu(x, q^{2}x)$. Let us note that if $y \in\{x, qx\}$ and $y\neq x$, then $qy \notin\{x, qx\}$. Consequently, we can say that the sectional curvature $\mu(x, qx)$ depends on $\varphi=\angle(x, qx)$. Analogously, $\mu(x, q^{2}x)$ depends on $\theta=\angle(x, q^{2}x)$.

We denote $\mu(x, qx)=\mu_{1}(\varphi)$ and $\mu(x, q^{2}x)=\mu_{2}(\theta)$.

\begin{theorem}\label{th4}
Let $(M, g, q)$ be a manifold with property \eqref{R}. If vectors $x$ and $u$ induce $q$-bases and $\{x, qx, q^{2}x, q^{3}x\}$ is orthonormal, then
\begin{equation}\label{mu-r2}
\begin{split}
   \mu_{1}(\varphi)-\mu_{1}(\frac{\pi}{2})&=\frac{\cos\varphi}{1-\cos^{2}\varphi}\Big(-2R(x, qx, q^{2}x, x)\\&+2(\cos\varphi) R(x, qx, qx, q^{2}x)\\&-(\cos\varphi) R(qx, q^{2}x, q^{3}x, x)\\&-2R(qx, q^{2}x, q^{2}x, x)\Big),
   \end{split}
\end{equation}
\begin{equation}\label{mu-r}
\begin{split}
   \mu_{2}(\theta)-\mu_{2}(\frac{\pi}{2})&=\frac{2\cos\theta}{1-\cos^{2}\theta}\Big(-2R(x, qx, qx, q^{2}x)\\&+(\cos\theta) R(qx, q^{2}x, q^{3}x, x)\Big),
   \end{split}
\end{equation}
where $\varphi=\angle(u, qu)$, $\theta=\angle(u, q^{2}u)$.
\end{theorem}
\begin{proof}
In \eqref{Rv1} we substitute

    1)\ $qx$ for $y$, $q^{2}x$ for $z$ and $x$ for $u$, then
    \begin{equation}\label{dop3}
    R(x, qx, q^{2}x, x)=R(q^{2}x, q^{3}x, x, q^{2}x)=R(q^{3}x, x, qx, q^{3}x),
    \end{equation}

    2)\  $qx$ for $y$, $qx$ for $z$ and $q^{2}x$ for $u$, and we have
    \begin{equation}\label{dop5}
    R(x, qx, qx, q^{2}x)=R(q^{2}x, q^{3}x, q^{3}x, x)=R(q^{3}x, x, x, qx),
    \end{equation}

    3)\  $qx$ for $y$, $q^{2}x$ for $z$ and $q^{3}x$ for $u$, then
    \begin{equation}\label{dop6}
    R(qx, q^{2}x, q^{3}x, x)=R(x, qx, q^{2}x, q^{3}x),
    \end{equation}

    4)\ $qx$ for $y$, $qx$ for $z$ and $q^{3}x$ for $u$, and we get
    \begin{equation}\label{dop7}
    R(qx, q^{2}x, q^{2}x, x)=R(x, qx, qx, q^{3}x)=R(q^{3}x, x, x, q^{2}x),
    \end{equation}

    5)\  $q^{2}x$ for $y$, $qx$ for $z$ and $q^{3}x$ for $u$, and we find
    \begin{equation}\label{dop8}
    R(x, q^{2}x, qx, q^{3}x)=0.
    \end{equation}

 Let $u=\alpha x+\beta qx +\gamma q^{2}x+\delta q^{3}x$, where $\alpha,\beta,\gamma, \delta \in \mathbb{R}$. From \eqref{f4} we get $qu=\delta x+\alpha qx +\beta q^{2}x + \gamma q^{3}x$, $q^{2}u=\gamma x+\delta qx +\alpha q^{2}x +\beta q^{3}x$ and \\$q^{3}u=\beta x+\gamma qx +\delta q^{2}x + \alpha q^{3}x$.
 Then, by using the linear properties of the curvature tensor $R$ and having in mind \eqref{Rv1}, \eqref{dop2}, \eqref{dop4}, \eqref{dop3} -- \eqref{dop8}, we obtain \begin{align*}
    R(u,qu,u,qu)&=\Big((\alpha^{2}-\beta\delta)^{2}+(\delta^{2}-\alpha\gamma)^{2}+(\beta^{2}-\alpha\gamma)^{2}+(\gamma^{2}-\beta\delta)^{2}\Big)R_{1}\\&+2\Big((\alpha\beta-\gamma\delta)(\gamma^{2}-\alpha^{2})+(\beta\gamma-\delta\alpha)(\delta^{2}-\beta^{2})\Big)R_{2}\\&+\Big((\alpha\beta-\gamma\delta)^{2}+(\beta\gamma-\delta\alpha)^{2}\Big)R_{3}\\&+2(\alpha^{2}+\gamma^{2}-2\beta\delta)(\delta^{2}+\beta^{2}-2\alpha\gamma)R_{4}\\&+2\Big((\alpha^{2}-\beta\delta)(\gamma^{2}-\beta\delta)+(\beta^{2}-\alpha\gamma)(\delta^{2}-\alpha\gamma)\Big)R_{5}\\&+2\Big((\alpha\beta-\gamma\delta)(\delta^{2}-\beta^{2})+(\beta\gamma-\delta\alpha)(\alpha^{2}-\gamma^{2})\Big)R_{6},\\ R(u,q^{2}u,u,q^{2}u)&=2\Big((\alpha\delta-\beta\gamma)^{2}+(\alpha\beta-\gamma\delta)^{2}\Big)R_{1}\\&+4\Big((\alpha\delta-\beta\gamma)(\gamma^{2}-\alpha^{2})+(\alpha\beta-\delta\gamma)(\delta^{2}-\beta^{2})\Big)R_{2}\\&+\Big((\alpha^{2}-\gamma^{2})^{2}+(\beta^{2}-\delta^{2})^{2}\Big)R_{3}\\&-2\Big((\alpha\beta-\gamma\delta)^{2}+(\beta\gamma-\alpha\delta)^{2}\Big)R_{5}\\&+4\Big((\beta^{2}-\delta^{2})(\alpha\delta-\beta\gamma)+(\alpha\beta-\gamma\delta)(\gamma^{2}-\alpha^{2})\Big)R_{6} ,
\end{align*}
where
\begin{equation}\label{r1-6}
    \begin{split}
R_{1}&=R(x, qx, x, qx), \qquad R_{2}=R(x, qx, q^{2}x, x),\\ R_{3}&=R(x, q^{2}x, x, q^{2}x),\quad
R_{4}=R(x, qx, qx, q^{2}x),\\ R_{5}&=R(qx, q^{2}x, q^{3}x, x),\ R_{6}=R(qx, q^{2}x, q^{2}x, x).
\end{split}
\end{equation}
Then
\begin{equation}\label{r+r}
    \begin{split}
   R(u,qu,u,qu)+\frac{1}{2}R(u,q^{2}u,u,q^{2}u)&=K_{1} R_{1}+ K_{2} R_{2}+K_{3} R_{3}\\&+K_{4} R_{4}+K_{5} R_{5}+K_{2} R_{6},
    \end{split}
\end{equation}
where \begin{align} \label{K}\nonumber K_{1}&=(\alpha^{2}-\beta\delta)^{2}+(\delta^{2}-\alpha\gamma)^{2}+(\beta^{2}-\alpha\gamma)^{2}+(\gamma^{2}-\beta\delta)^{2}\\\nonumber &+(\alpha\delta-\beta\gamma)^{2}+(\alpha\beta-\gamma\delta)^{2},\\ \nonumber K_{2}&=2\Big((\alpha\beta-\gamma\delta)(\gamma^{2}-\alpha^{2})+(\beta\gamma-\delta\alpha)(\delta^{2}-\beta^{2})\\\nonumber &+(\alpha\delta-\beta\gamma)(\gamma^{2}-\alpha^{2})+(\alpha\beta-\delta\gamma)(\delta^{2}-\beta^{2})\Big),\\ K_{3}&=(\alpha\beta-\gamma\delta)^{2}+(\beta\gamma-\delta\alpha)^{2}\\\nonumber &+\frac{1}{2}\Big((\alpha^{2}-\gamma^{2})^{2}+(\beta^{2}-\delta^{2})^{2}\Big), \end{align}\begin{align*}
K_{4}&=2(\alpha^{2}+\gamma^{2}-2\beta\delta)(\delta^{2}+\beta^{2}-2\alpha\gamma),\\ K_{5}&=2\Big((\alpha^{2}-\beta\delta)(\gamma^{2}-\beta\delta)+(\beta^{2}-\alpha\gamma)(\delta^{2}-\alpha\gamma)\Big)\\ &-(\alpha\beta-\gamma\delta)^{2}-(\beta\gamma-\alpha\delta)^{2}.\end{align*}
Since the $q$-basis $\{x, qx, q^{2}x, q^{3}x\}$ is orthonormal, we have
\begin{equation}\nonumber
    g(u, u)= g(qu, qu)=\alpha^{2}+\beta^{2}+\gamma^{2}+\delta^{2},
     \end{equation}
     \begin{equation}\nonumber
g(u,qu)= \alpha\delta+\alpha\beta+\beta\gamma+\delta\gamma,\quad g(u, q^{2}u)=2(\alpha\gamma+\delta\beta).
\end{equation}
Due to \eqref{2.1} and \eqref{3.3} we get
\begin{equation}\nonumber
    \mu(u,qu)=\frac{R(u, qu, u, qu)}{g^{2}(u, u)-g^{2}(u, qu)}\ ,\quad\mu(u,q^{2}u)=\frac{R(u, q^{2}u, u, q^{2}u)}{g^{2}(u, u)-g^{2}(u, q^{2}u)}\ .
\end{equation}
We suppose that $g(u, u)=1$ and we obtain
     \begin{equation}\label{mu3}
       \mu(u,qu)=\frac{R(u, qu, u, qu)}{1-\cos^{2}\varphi}\ ,\quad \mu(u,q^{2}u)=\frac{R(u, q^{2}u, u,  q^{2}u)}{1-\cos^{2}\theta}\ ,
 \end{equation}
\begin{equation}\label{alfa-delta}
   \alpha^{2}+\beta^{2}+\gamma^{2}+\delta^{2}=1,\
   \alpha\delta+\alpha\beta+\beta\gamma+\delta\gamma=\cos\varphi, \ 2(\alpha\gamma+\delta\beta)=\cos\theta.
\end{equation}
 From \eqref{alfa-delta} we express $\alpha, \beta, \gamma, \delta$ by $\cos\varphi$ and $\cos\theta$. Then, taking into account \eqref{K}, we get
    \begin{align*}
    % \nonumber to remove numbering (before each equation)
     &K_{1}=1-\cos^{2}\varphi,\quad K_{2}=-2\cos\varphi(1-\cos\theta),\quad K_{3}=\frac{1}{2}(1-\cos^{2}\theta), \\& K_{4}=2(-\cos\theta+\cos^{2}\varphi), \quad K_{5}=\cos^{2}\theta
-\cos^{2}\varphi.
    \end{align*}
Thus, \eqref{r+r} and \eqref{mu3} imply
\begin{equation}\label{mu+mu}
    \begin{split}
   (1-\cos^{2}\varphi)\mu(u,qu)&+\frac{1}{2}(1-\cos^{2}\theta)\mu(u,q^{2}u)=\\&(1-\cos^{2}\varphi)R_{1} -2\cos\varphi(1-\cos\theta)R_{2}\\&+\frac{1}{2}(1-\cos^{2}\theta)R_{3}+2(-\cos\theta+\cos^{2}\varphi)R_{4}\\&+(\cos^{2}\theta
-\cos^{2}\varphi)R_{5}-2\cos\varphi(1-\cos\theta)R_{6}.
    \end{split}
\end{equation}
In \eqref{mu+mu} first we substitute $\varphi=\frac{\pi}{2}$ and then $\theta=\frac{\pi}{2}$. Thus we obtain \eqref{mu-r} and \eqref{mu-r2}.
   \end{proof}
   \begin{theorem}\label{th6}
Let $(M,g, q)$ be a manifold with property \eqref{R}. If a vector $u$ induces a $q$-basis, then the following equality is valid
\begin{equation}\label{mu-r4}
\begin{split}
   \mu_{1}(\varphi)&=\frac{1}{1-\cos^{2}\varphi}\Big((1-4\cos^{2}\varphi)\mu_{1}(\frac{\pi}{2})\\&+\frac{3}{4}(\cos\varphi+2\cos^{2}\varphi)\mu_{1}(\frac{\pi}{3})\\&+\frac{3}{4}(2\cos^{2}\varphi-\cos\varphi)\mu_{1}(\frac{2\pi}{3})\Big),
   \end{split}
\end{equation}
where  $\varphi=\angle(u, qu)$.
\end{theorem}
\begin{proof}
In \eqref{mu-r2} first we substitute $\varphi=\frac{\pi}{3}$ and then $\varphi=\frac{2\pi}{3}$. Due to \eqref{r1-6} and having in mind that $\{x, qx, q^{2}x, q^{3}x\}$ is an orthonormal $q$-basis, we get
\begin{align*}
&3\Big(\mu_{1}(\frac{\pi}{3})-\mu_{1}(\frac{\pi}{2})\Big)=4\Big(-R_{2}+\frac{1}{2}R_{4}-\frac{1}{4}R_{5}-R_{6}\Big),\\
&3\Big(\mu_{1}(\frac{2\pi}{3})-\mu_{1}(\frac{\pi}{2})\Big)=4\Big(R_{2}+\frac{1}{2}R_{4}-\frac{1}{4}R_{5}+R_{6}\Big).
\end{align*}
From the latter equalities we find the tensors $R_{4}-\frac{1}{2}R_{5}$ and $R_{2} + R_{6}$. Then \eqref{mu-r2} implies \eqref{mu-r4}.
 \end{proof}
\begin{theorem}\label{th5}
Let $(M,g, q)$ be a manifold with property \eqref{R}. If a vector $u$ induces a $q$-basis, then the following equality is valid
\begin{equation}\label{mu-r3}
\begin{split}
   \mu_{2}(\theta)&=\frac{1}{1-\cos^{2}\theta}\Big((1-4\cos^{2}\theta)\mu_{2}(\frac{\pi}{2})\\&+\frac{3}{4}(\cos\theta+2\cos^{2}\theta)\mu_{2}(\frac{\pi}{3})\\&+\frac{3}{4}(2\cos^{2}\theta-\cos\theta)\mu_{2}(\dfrac{2\pi}{3})\Big),
   \end{split}
\end{equation}
where $\theta=\angle(u, q^{2}u)$.
\end{theorem}
\begin{proof}
In \eqref{mu-r} first we substitute $\theta =\frac{\pi}{3}$ and then $\theta =\frac{2\pi}{3}$. Thus we get
\begin{align*}
&3\Big(\mu_{2}(\frac{\pi}{3})-\mu_{2}(\frac{\pi}{2})\Big)=2R_{5}-8R_{4} ,\\
&3\Big(\mu_{2}(\frac{2\pi}{3})-\mu_{2}(\frac{\pi}{2})\Big)=2R_{5}+8R_{4} .
\end{align*}
  Taking into account the last system and \eqref{mu-r}, we obtain \eqref{mu-r3}.
 \end{proof}
 \begin{theorem}\label{th7}
Let $(M,g, q)$ be a manifold with property \eqref{R1}. If a vector $u$ induces a $q$-basis, then the following equalities are valid
\begin{equation}\label{mu-r5}
   \mu_{2}(\theta)=0,\quad \mu_{1}(\varphi)=\frac{(1-\cos\theta)^{2}}{1-\cos^{2}\varphi}\mu_{1}(\frac{\pi}{2}),
\end{equation}
where $\theta=\angle(u, q^{2}u)$ and $\varphi=\angle(u, qu)$.
\end{theorem}
\begin{proof}
From \eqref{R1} we have
\begin{equation}\label{q-q2-q3}
  R(x, y, qz, qu)=R(x, y, q^{2}z, q^{2}u).
\end{equation}

In \eqref{R1}, \eqref{q-q2-q3} we substitute

1) $qx$ for $y$, $x$ for $z$, $qx$ for $u$, and we get
 \begin{equation}\label{r11}
 R(x,qx,x,qx)=R(x,qx,qx,q^{2}x )=R(x,qx,q^{2}x,q^{3}x),
 \end{equation}

2) $q^{2}x$ for $y$, $x$ for $z$, $q^{2}x$ for $u$, and we obtain
  \begin{equation}\label{r12}
    R(x, q^{2}x, x, q^{2}x)=R(x, q^{2}x, qx, q^{3}x)=R(x, q^{2}x, q^{2}x, x),
    \end{equation}

3) $qx$ for $y$, $q^{2}x$ for $z$, $x$ for $u$, then
 \begin{equation}\label{r13}
 R(x, qx,q^{2}x,x)=R(x,qx,q^{3}x,qx)=R(x,qx,x,q^{2}x).
 \end{equation}
Comparing the identities \eqref{dop6}, \eqref{dop7}, \eqref{r11}, \eqref{r12}, \eqref{r13} and having in mind \eqref{r1-6}, we obtain $R_{1}=R_{4}=R_{5}$ and $R_{2}=R_{3}=R_{6}=0$.
 Then \eqref{mu+mu} implies \eqref{mu-r5}.
\end{proof}
We note that Proposition~4.2 from \cite{1} is a particular case of Theorem~\ref{th7}.

\section{A Lie group as a $4$-dimensional Riemannian manifold with a circulant structure}\label{sec:3}

Let $G$ be a $4$-dimensional real connected Lie group and $\mathfrak{g}$ be its Lie algebra with a basis $\{x_{1}, x_{2},x_{3},x_{4}\}$. We introduce a structure $q$ and left invariant metric $g$ as follows
\begin{equation}\label{lie}
  qx_{1}=x_{2} ,\ qx_{2}=x_{3},\ qx_{3}=x_{4},\ qx_{4}=x_{1},
\end{equation}
  \begin{equation}\label{g}
 g(x_{i}, x_{j})= \left\{ \begin{array}{ll}
                        0, & i\neq j \hbox{;} \\
                        1, & i=j \hbox{.}
                      \end{array}
                    \right.
  \end{equation}

Obviously, \eqref{f4} and \eqref{2.1} are valid.
Therefore $(G, g, q)$ is a Riemannian manifold with \eqref{f4} and \eqref{2.1}.

For the manifold $(G, g, q)$ we suppose that $g$ is a Killing metric, i.e.
\begin{equation}\label{killing}
    g([x_{i}, x_{j}],x_{k})+ g([x_{i}, x_{k}],x_{j})=0.
\end{equation}
According to \eqref{killing} and the Jacobi identity for the commutators $[x_{i}, x_{j}]$ we obtain
\begin{align}\label{skobki2}\nonumber
  [x_{1}, x_{2}]&=\lambda_{1}x_{3}+\lambda_{2}x_{4},\qquad [x_{1}, x_{3}]=-\lambda_{1}x_{2}+\lambda_{4}x_{4},\\ [x_{2}, x_{3}]&=\lambda_{1}x_{1}+\lambda_{3}x_{4} ,\qquad
[x_{1}, x_{4}]=-\lambda_{2}x_{2}-\lambda_{4}x_{3},\\\nonumber [x_{2}, x_{4}]&=\lambda_{2}x_{1}-\lambda_{3}x_{3},\qquad [x_{3}, x_{4}]=\lambda_{4}x_{1}+\lambda_{3}x_{2},
\end{align}
where $\lambda_{i}\in \mathbb{R}$.

Vice versa, if \eqref{skobki2} are valid for a Riemannian manifold $(G, g, q)$, where the structure $q$ and the metric $g$ on the Lie group $G$ are determined by \eqref{f4} and \eqref{2.1}, then the Jacobi identity for commutators $[x_{i}, x_{j}]$ is satisfied and the metric $g$ is Killing.

Therefore, we establish the truthfulness of the following
\begin{theorem}\label{kt}
Let $(G, g, q)$ be a $4$-dimensional Riemannian manifold, where $G$ is the connected Lie group with an associated Lie algebra $\mathfrak{g}$, determined by a global basis $\{x_{i}\}$ of left invariant vector fields, and $q$ and $g$ are the structure and the metric, determined by \eqref{lie} and \eqref{g}. Then $(G, g, q)$ is a Riemannian manifold with a circulant structure $q$ and a Killing metric $g$, which satisfy \eqref{f4} and \eqref{2.1} if and only if $G$ belongs to a Lie group, determined by \eqref{skobki2}.
\end{theorem}

Further, $(G, g, q)$ will stand for the Riemannian manifold determined by the conditions of Theorem~\ref{kt}.

Since $g$ is a Killing metric, then the components of $R$ are (\cite{32})

\begin{equation}\label{r}
    R_{ijkh}=-\frac{1}{4}g\Big([x_{i},x_{j}],[x_{k},x_{h}]\Big).
\end{equation}

According to \eqref{g}, \eqref{skobki2} and \eqref{r} we calculate the following components of the curvature tensor $R$:
\begin{align}\label{r1}\nonumber
    R_{1212}&=-\frac{1}{4}(\lambda_{1}^{2}+\lambda_{2}^{2}),\quad R_{1414}=-\frac{1}{4}(\lambda_{2}^{2}+\lambda_{4}^{2}),\\\nonumber
R_{2323}&=-\frac{1}{4}(\lambda_{1}^{2}+\lambda_{3}^{2}), \quad R_{3434}=-\frac{1}{4}(\lambda_{3}^{2}+\lambda_{4}^{2}),\\
R_{1313}&=-\frac{1}{4}(\lambda_{1}^{2}+\lambda_{4}^{2}),\quad R_{2424}=-\frac{1}{4}(\lambda_{2}^{2}+\lambda_{3}^{2}),\\\nonumber
R_{1213}&=R_{2434}=-\frac{1}{4}\lambda_{2}\lambda_{4}, \quad R_{2324}=R_{1314}=-\frac{1}{4}\lambda_{1}\lambda_{2},\\\nonumber
R_{1424}&=R_{1323}=-\frac{1}{4}\lambda_{3}\lambda_{4}, \quad R_{3134}=R_{2124}=-\frac{1}{4}\lambda_{1}\lambda_{3},\\\nonumber
R_{1214}&=R_{3234}=\frac{1}{4}\lambda_{1}\lambda_{4}, \qquad R_{1434}=R_{2123}=\frac{1}{4}\lambda_{2}\lambda_{3}.
\end{align}
The rest of nonzero components are obtained from the properties
$$R_{ijks}=R_{ksij},\ R_{ijks}=-R_{jiks}=-R_{ijsk}.$$

\begin{prop}\label{kt2}
Let $(G, g, q)$ be a manifold determined by the conditions of Theorem~\ref{kt}.
  Then $(G, g, q)$  satisfies the identity \eqref{R} if and only if
\begin{equation}\label{landa}
\lambda_{1}=\varepsilon\lambda_{2}=\varepsilon\lambda_{3}=\lambda_{4},\ \varepsilon=\pm 1.
\end{equation}
\end{prop}
\begin{proof}

According to \eqref{lie} we obtain that \eqref{R} is equivalent to
\begin{align*}\nonumber
  R_{1212}=R_{3434}=R_{2323}=R_{1414},\quad R_{1313}=R_{2424},\\
R_{1213}=R_{2324}=R_{1424}=R_{3134},\
R_{1214}=R_{1434}=R_{2123}=R_{3234}, \\\nonumber
R_{1224}=R_{3123}=R_{3114}=R_{4234}, \ R_{1324}=0.
\end{align*}
Then \eqref{r1} implies
\begin{align*}
    \lambda_{2}=\lambda_{3},\ \lambda_{1}=\lambda_{4},\ \lambda_{1}^{2}+\lambda_{4}^{2}=\lambda_{2}^{2}+\lambda_{3}^{2},\ \lambda_{1}\lambda_{4} = \lambda_{2}\lambda_{3}.
\end{align*}
So we obtain \eqref{landa}.

\end{proof}
From \eqref{r1} and \eqref{landa} we calculate
\begin{align}\label{rlamda}\nonumber
  R_{1212}=R_{1414}=R_{2323}&=R_{3434}=R_{1313}=R_{2424}=-\frac{1}{2}\lambda_{1}^{2},\\
R_{1213}=R_{2434}=R_{2324}&=R_{1314}=R_{1424}=\\\nonumber &R_{3431}=R_{1323}=R_{2124}=-\frac{1}{4}\varepsilon\lambda_{1}^{2}, \\\nonumber
R_{1214}=R_{3234}=R_{1434}&=R_{2123}=\frac{1}{4}\lambda_{1}^{2}.
\end{align}
Having in mind \eqref{rlamda} and the formulas
\begin{equation*}
    \rho(y,z)=g^{ij}R(e_{i}, y, z, e_{j}),\qquad
    \tau=g^{ij}\rho(e_{i}, e_{j}),
\end{equation*}
we get the components of the Ricci tensor $\rho$ and the value of the scalar curvature $\tau$ as follows:
\begin{align}\label{rho}\nonumber
    \rho_{11}=\rho_{22}=\rho_{33}=\rho_{44}=\frac{3}{2}\lambda_{1}^{2},\\ \rho_{12}=\rho_{14}=\rho_{23}=\rho_{34}=\frac{1}{2}\varepsilon\lambda_{1}^{2}, \\\nonumber
\rho_{13}=\rho_{24}=-\frac{1}{2}\lambda_{1}^{2},
\end{align}
\begin{equation}\label{tau2}
 \tau=6\lambda_{1}^{2}.
\end{equation}

By using \eqref{3.3} for the sectional curvatures of the basic $2$-planes we find
\begin{equation}\label{mu4}
    \mu(x_{1},x_{2})=\mu(x_{1},x_{4})=\mu(x_{2},x_{3})= \mu(x_{2},x_{4})=\mu(x_{1},x_{3})=-\frac{1}{2}\lambda_{1}^{2}.
\end{equation}

Therefore, we arrive at the following

\begin{prop}\label{kt3}
Let $(G, g, q)$ be a manifold determined by the conditions of Theorem~\ref{kt}.
 If $(G, g, q)$  satisfies the identity \eqref{R}, then
\begin{itemize}
\item[(i)]  The nonzero components of the curvature tensor $R$ and the Ricci tensor $\rho$ are \eqref{rlamda} and \eqref{rho};

 \item[(ii)]  The scalar curvature $\tau$ is \eqref{tau2};

\item[(iii)]  $(G, g, q)$ is of constant sectional curvatures \eqref{mu4}, i.e. $(G, g, q)$  is conformally flat manifold.
\end{itemize}
\end{prop}
According to \eqref{lie} we obtain that \eqref{R1} is equivalent to the equalities
\begin{align*}
  R_{1212}=R_{1414}=R_{2323}=R_{3434}&=R_{1223}=R_{1241}=\\R_{4134}&=R_{1234}=R_{2334}= R_{2341},\\
R_{1313}=R_{2424}=R_{1324}=R_{1213}&=R_{1224}=R_{1431}=\\
R_{2441}&=R_{2423}=R_{2331}=R_{1334}=R_{2434}=0.
\end{align*}
Then, by using \eqref{r1} and \eqref{tau2} we have the following
\begin{prop}
  Let $(G, g, q)$ be a manifold determined by the conditions of Theorem~\ref{kt}. Then the following propositions are equivalent:
\begin{itemize}
  \item[(i)] $\lambda_{1}=\lambda_{2}=\lambda_{3}=\lambda_{4}=0$, i.e. $G$ is abelian;
  \item[(ii)] $(G, g, q)$  satisfies the identity \eqref{R1};
  \item[(iii)] $\tau=0$, i.e. $(G, g, q)$  is a scalar flat manifold with respect to $\nabla$.
\end{itemize}

\end{prop}

%------------------------------------------------------------------
\section*{Acknowledgments}

This work was partially supported by project NI15-FMI-004 of the
Scientific Research Fund, Paisii Hilendarski University of
Plovdiv, Bulgaria.


\begin{thebibliography}{99}


\bibitem{8}
\textsc{Davis, P. J.}: Circulant matrices,  New York, John Wiley and Sons, 250 (1979)

\bibitem{22}
   \textsc{Gadea, P. M., Montesinos, A.}: Spaces of constant para-holomorphic sectional curvature, Pacific J. Math., \textbf{136}, No~1, 85--101 (1989)

   \bibitem{132}
\textsc{Gray, A., Vanhecke, L.}: Almost Hermitian manifolds with constant holomorphic sectional curvature, Appl. Math., \textbf{104}, 170--179 (1979)

   \bibitem{152}
\textsc{Gray, R. M.}: Toeplitz and circulant matrices: A review, Found. Trends Commun. Inf. Theory, \textbf{2}, No~3, 155--239 (2006)

   \bibitem{32}
\textsc{Mekerov, D.}: Lie groups as $4$-dimensional Riemannian or pseudo-Riemann\-ian almost product manifolds with nonitegrable structure, J. Geom.,  \textbf{90},  165--174 (2008)

\bibitem{13}
\textsc{Muzychuk, M.}: A solution of the isomorphism problem for circulant graphs, In: Proc. London Math. Soc., \textbf{88}, No~3, 1--41 (2004)

\bibitem{72}
   \textsc{Olson, B., Shaw, S., Shi, C., Pierre C., Parker, R. G.}: Circulant matrices and their application to vibration analysis, Appl. Mech. Rev., \textbf{66}, No~4, 1--41 (2014)

\bibitem{1}
 \textsc{Razpopov, D.}: Four-dimensional Riemannian manifolds with two circulant structures, In: Proc. of $44$-th Spring Conf. of UBM, SOK "Kamchia", Bulgaria, 179--185 (2015)

 \bibitem{2}
 \textsc{Roth, R. M., Lempel, A.}: Application of circulant matrices to the construction and decoding of linear codes, IEEE Trans. Inf. Theory, \textbf{36}, No~5, 1157--1163 (1990)

\bibitem{52}
\textsc{Staikova, M., Gribachev, K., Mekerov, D.}: Riemannian $P$-manifolds of constant sectional curvatures, Serdica Math. J., \textbf{17}, 212--219 (1991)

\bibitem{51}
\textsc{Staikova, M., Gribachev, K.}: Canonical conections and their conformal invariants on Riemannian $P$-manifolds, Serdica Math. J., \textbf{18}, 150--161 (1992)

\bibitem{11}
\textsc{Yano, K.}: Differential geometry on complex and almost complex spaces. Pure and Applied Math. \textbf{49},  New York, Pergamont Press Book, 326 (1965)


% and use \bibitem to create references. Consult the Instructions
% for authors for reference list style.
%
\end{thebibliography}
\end{document}